\documentclass{article}
\usepackage{amsmath,hyperref}
\usepackage{amssymb}
\usepackage{amsthm}
\usepackage{romannum}
\usepackage{graphicx}
\usepackage[dvipsnames]{xcolor}
\usepackage{comment}
\usepackage{hyperref}
  \newtheorem{theorem}{Theorem}[section]
        \newtheorem{proposition}[theorem]{Proposition}
        \newtheorem{lemma}[theorem]{Lemma}
        \newtheorem{corollary}[theorem]{Corollary} 
         
        \newtheorem{claim}{Claim} 
        \newtheorem{remark}[theorem]{Remark} 
\numberwithin{equation}{section}

\title{Global Existence for Coupled  3-D  Nonlinear Wave and Klein-
Gordon Equations with Large Derivatives of Initial Data}
\author{Guocong Shang\footnote{Southern University of Science and Technology, Department of Mathematics. \quad Email: 12432027@mail.sustech.edu.cn}} 
\date{\today}
\begin{document}
\maketitle
\begin{abstract}
    We consider the Cauchy problem of coupled 3-D wave and Klein-Gordon equations with a quadratic form of nonlinearity. We show global existence under several conditions, including large derivative data for wave equations and the null conditions.
\end{abstract}
\section{Introduction}
We consider the following coupled wave and Klein-Gordon equations
\begin{equation}
\left\{
\begin{aligned}
&-\Box u=Q_0(u,w), \\
&-\Box w +  w = Q_0(u,w)
\end{aligned} \label{add 1.2}
\right.
\end{equation}
in $\mathbb{R}_+ \times \mathbb{R}^3=\{(t,x)|t>0, x\in \mathbb{R}^3\}$. The flat box operator $-\Box:=\partial_t^2-\sum_{i=1}^3\partial_{x_i}^2$, the nonlinearities $Q_0(u, w) := \partial_\alpha u \, \partial^\alpha w = -\partial_t u \, \partial_t w + \sum_{i=1}^3 \partial_{x_i} u \, \partial_{x_i} w$.

The initial conditions are given at $t=0$:
\begin{equation*}
\begin{split}
    &u(0, x) = u_0(x), \quad (\partial_t u)(0, x) = u_1(x),\quad x \in \mathbb{R}^3,
    \\ &  \label{add.1.2}
w(0, x) = w_0(x), \quad (\partial_t w)(0, x) = w_1(x), \quad x \in \mathbb{R}^3.
\end{split} 
\end{equation*}
We will show the system \eqref{add 1.2} admits global solutions with large derivatives initial data of $u$ in the $L^2$ sense, precisely, 
\begin{theorem}
    Let $N>40$ be a fixed integer, $0<p<\frac{1}{10}$, there exists an $\epsilon_0>0$ such that $\forall 0<\epsilon<\epsilon_0$ for all  initial data satisfying the following conditions, 
    \begin{equation}
   \begin{split}
     & \sum_{|I|\le N} ||\langle x\rangle^{|I|}\nabla^I u_0||_{L^2}\leq \frac{100}{\epsilon},
    \\  & \sum_{|I|\le N} ||\langle x\rangle^{|I|+1}\nabla^I u_1||_{L^2}+||\langle x\rangle^{|I|+1}\nabla\nabla^I u_0||_{L^2}\leq 1,
       \\&  \sum_{|I|\le N} ||\langle x\rangle^{|I|}\nabla\nabla^I u_1||_{L^2}+   ||\langle x\rangle^{|I|}\nabla \nabla\nabla^I u_0||_{L^2}\leq \epsilon,
       \\& \sum_{|I|\le N} ||\langle x\rangle^{|I|+8}\nabla^I w_0||_{L^2}\leq \epsilon,
       \\  & \sum_{|I|\le N} ||\langle x\rangle^{|I|+7}\nabla\nabla^I w_0||_{L^2}+ ||\langle x\rangle^{|I|+7}\nabla^I w_1||_{L^2}\leq \epsilon.\label{1.3}
   \end{split} 
    \end{equation}
   Then Cauchy problem associates to \eqref{add 1.2}, \eqref{1.3} admits a global-in-time solution with asymptotic behavior
   \begin{equation}
   \begin{split}
      & %|u(t,x)|\lesssim  \langle t+r \rangle ^{-1+p},
      |\partial u(t,x)|\lesssim   \langle t+r \rangle ^{-1} \langle t-r \rangle ^{-\frac{1}{8}},  
      \\ &  |w(t,x)|\lesssim  \epsilon \langle t+r \rangle ^{-\frac{3}{2}},  |\partial w(t,x)|\lesssim   \epsilon \langle t+r \rangle ^{-\frac{3}{2}}. \label{thm}
   \end{split}
   \end{equation}

\end{theorem}
\begin{remark}
    To the best of the author's knowledge, this setting for the initial data has appeared in \cite{dong2023cubic}.
\end{remark}
\begin{remark}
    The requirements for $u_0$ and $u_1$ in \eqref{1.3} are to guarantee the existence of such initial data. The Theorem 1.1 shows it's admissible to require only boundedness of $\|\partial u(0,\cdot)\|_{L^2}$ for the global existence of \eqref{add 1.2}, while the smallness for $\|\partial \partial u(0,\cdot)\|_{L^2}$ is needed. The requirements for $w_0$ and $w_1$ are the usual smallness requirements.
\end{remark}

\subsection{Relevant result revisited} 

 The coupled
 wave and Klein-Gordon systems can be derived from many important physical models,
 like the Dirac-Proca model and Klein-Gordon-Zakharov model, the Maxwell-Klein-Gordon
 model, the Einstein-Klein-Gordon mode, and many others. We first go through some relevant results of the nonlinear wave equation in $[0,\infty)\times \mathbb{R}^n (n\ge 2)$:
 \begin{equation}
     \Box u = g^{ij} \partial_k u \partial_i \partial_j u, \label{wave}
 \end{equation}
  where $g^{ij}$ are constants and the initial data are sufficiently nice. When $n\ge4$, the GWP of \eqref{wave} is obtained in \cite{hormander1997lectures}. When $n=3$, the GWP of \eqref{wave} is obtained in the celebrated work of Klainerman \cite{klainerman1986null} and Christodoulou \cite{christodoulou1986global} under the null condition $g^{kij} w_k w_i w_j = 0,\text{for   }  w_0=-1, (w_1,w_2,w_3)\in S^2$. If one removes the requirements of the null condition, \eqref{wave} might blow up in finite time (cf.\cite{john1979blow}). When $n=2$, the slow decay rate in time causes the major difficulty. Alinhac introduced the "ghost weight" method in \cite{alinhac2000null} and
   \cite{alinhac2001null} to show the \eqref{wave} has global-in- time solution with small initial data. Additionally, the highest norm of the solution exhibits polynomial growth in time.

The nonlinear Klein-Gordon equations have also been extensively studied. The work by Klainerman \cite{klainerman1987global} and Shatah \cite{shatah1985normal} establishes the global existence of solutions to the nonlinear Klein-Gordon equations. After that, various exciting results on nonlinear wave equations, the nonlinear Klein
Gordon equations, and their coupled systems, come out. Bachelot \cite{bachelot1988probleme} considered the Dirac-Klein-Gordon system.
 Later on, 
 Georgiev \cite{georgiev1990global} proved the global existence result for the nonlinearities satisfying the strong
 null condition (i.e. $\partial_\alpha u\partial_\beta v - \partial_\beta u\partial_\alpha v, \alpha,\beta \in \{0,1,2,3\}$ ), which are compatible with the
 vector fields, excluding the scaling vector field.

 A long time after that,
 Katayama \cite{katayama2012global,katayama2018global} proved a large class of coupled wave and Klein-Gordon equations. In
 \cite{katayama2012global}, for example, Katayama considered the bad nonlinearities of the type $uv,u\partial v$, but with the assumption that the nonlinearities in the wave equation of $u$ have divergence form, and this means that the bad nonlinearities are essentially of the type $\partial uv,\partial u\partial v$, which are not very bad.
Later on, LeFloch and Ma also tackled a large class of the coupled wave and Klein-Gordon equations with compactly supported initial data in \cite{lefloch2014hyperboloidal} using the hyperboloidal foliation method, which is essentially Klainerman’s vector field method on hyperboloids. 
 Recently, Dong \cite{dong2021global} showed that the coupled wave-Klein-Gordon system admits a global solution in 2D with a quadratic form $Q_{\alpha\beta}$ and small initial data. For
other recent results on wave-Klein-Gordon systems, see \cite{andersson2018wave}, \cite{boussaid2017spectral}, \cite{cacciafesta2024strichartz}, \cite{dai2025maxwell}, \cite{dong2023asymptotic}, \cite{ma2025non}, \cite{dong2020stability}, \cite{ma2021global}.

  \subsection{The difficulties and corresponding strategy} 
     Compared with the previous result, the result in this paper considers data with large derivatives(but only for the wave component) and nonlinearity $Q_0$, therefore it brings three parts of main difficulty: 1) when applying Klainerman's vector fields method, the lack of scaling vector field $S=t \partial_t + r \partial_r$ makes the decay rate of Klein-Gordon part much slower than needed. 2)To control the quadratic form $Q_0$, the vector field $S$ is needed(see Lemma 2.2) while quadratic form $Q_{\alpha\beta}$ does not, see \eqref{add 2.2}. 3) The largeness brought by the initial data of $u$ makes the bootstrap assumptions \eqref{3.1} hard to close.

To address the difficulties mentioned above, we use the Sobolev embedding $\iota: \dot H^1(\mathbb{R}^3)\rightarrow L^6(\mathbb{R}^3)$ and point-wise interpolation technique to balance largeness of the $\| u\|_{L^2}$ and $\|\partial u\|_{L^2}$ term and sufficient decay rate (see Proposition 3.3). The point-wise interpolation (see Proposition 3.4) is inspired by \cite{dong2023cubic}. They introduce this balancing idea to treat the Dirac equation. For the scaling vector field, we control $|Su|\le \langle t+r\rangle|\partial u|$ by brute force. However, this brute estimate introduces an extra $\langle t+r\rangle$ growth; the natural decay rate of the Klein-Gordon components can cancel it.

\subsection{Organization of the paper}
The paper is organized as follows. In section \ref{sec 2}, we introduce the notation and basic estimates that will be used in sections \ref{sec3}. The global existence of \eqref{add 1.2} is proved by a bootstrap argument. The set of bootstrap bounds, which are the energy bounds on wave and Klein-Gordon components, is given in \ref{sec3.1}. It will be improved in different steps in \ref{sec 3.4} respectively.

\section {Preliminaries} \label{sec 2}
\subsection{Basic notations}

In the $(1+3)$ dimensional spacetime, we adopt the signature $(-, +, +,+)$. We denote a point in $\mathbb{R}^{1+3}$ by $(x_0, x_1, x_2,x_3) = (t, x_1, x_2,x_3)$, and denote its spatial radius by $r = \sqrt{x_1^2 +x_2^2+ x_3^2}$. 

The Japanese bracket is written as \( \langle \cdot \rangle := \sqrt{1 + |\cdot|^2} \). Throughout the paper, we use \( A \lesssim B \) to indicate \( A \leq CB \) with \( C \) a universal constant, and \( A \sim B \) if both \( A \lesssim B \) and \( B \lesssim A \) hold. Through this paper all the $L^p$ , $1\le p\le \infty$  norms are taken in $\mathbb{R}^3$. The $\partial$ refers to spatial and time derivatives, and $\nabla$ refers to spatial derivatives only.

To apply Klainerman's vector fields method, we first introduce the vector fields:
\begin{itemize}
\item Translations: $\partial_\alpha$, \quad $\alpha = 0, 1, 2,3$.

\item Rotations: $\Omega_{ab} = x_a \partial_b - x_b \partial_a$,  \quad $a, b = 1, 2,3$.

\item Lorentz boosts: $L_a = x_a \partial_t + t \partial_a$, \quad $a = 1, 2,3$.

\item Scaling vector field: $ S = t \partial_t + r \partial_r=t \partial_t + \sum_{i=1}^3 \partial_{x_i}$.

\end{itemize}
We will use $Z$ to denote a general vector field (without the scaling vector field $S$) in 
$$
V := \{ \partial_\alpha, \Omega_{ab}, L_a \}.
$$
We also define vector fields, including scaling vector fields. 
$$ \Gamma:=\{S,L_a,\Omega_{ab},\partial_a\}.
$$
In addition, we also need good derivatives,
$$
G_a:=r^{-1}(x_a\partial _t+r\partial_a).
$$
Let $\left\{ \zeta_k \right\}_{k=0}^{\infty}$ be a Littlewood-Paley partition of unity, i.e.

\[
1 = \sum_{k=0}^{\infty} \zeta_k(\tau), \quad \tau \geq 0, \quad \zeta_k \in C_0^{\infty}(\mathbb{R}), \quad \zeta_k \geq 0 \text{ for all } k \geq 0,
\]
as well as

\[
\text{supp } \zeta_0 \cap [0, \infty) = [0, 2], \quad \text{supp } \zeta_k \subset [2^{k-1}, 2^{k+1}] \text{ for all } k \geq 1.
\]
\subsection{Estimates for commutators and null forms}
The following three lemmas for commutators will be frequently used, see \cite{sogge1995lectures, hou2020global}.
\begin{lemma} \label{lemma 2.1}
    For any $Z', Z'' \in V$ we have
\begin{equation}
\begin{split}
   & [\Box, Z'] = 0,
\quad
\big| [Z', Z''] w \big| \lesssim  \big| Z w \big|, \quad [-\Box,S]=-2\Box ,
\\ &
\big| [Z, \partial] w \big| + \big| [S, Z] w \big| \lesssim  \big| \partial w \big| ,\quad Z(fg)=Z(f)g+Z(g)f.
\label{2.1}
\end{split}
\end{equation}
In addition, if we act $\Gamma$  on null forms, we have
\begin{equation}
    \Gamma Q_0(u,v)-Q_0(\Gamma u,v)-Q_0(u,\Gamma v) =0.
\end{equation}
where $w,u,v$ are a sufficiently nice functions.%,  where $Q_{ab}(u,v)=\partial_a u\partial_bv-\partial_av\partial_bu$, $a,b\in\{1,2,3\}$.
\end{lemma}
\begin{lemma}
  For  sufficiently regular functions $u$, $v$, we have 
    \begin{equation}
    \begin{split}
        & |Q_{ab}(u,v)|+|Q_0(u,v)|\lesssim \sum_a \left( |G_a u| |\partial v| + |G_a v| |\partial u|\right),
         \\ & |Q_0(u,v)| \lesssim \langle t + |x|\rangle ^{-1} \left( |SuZ v| + |Z uZ v|\right),
          \\ & |Q_{ab}(u,v)| \lesssim \langle t + |x|\rangle ^{-1} 
          \left(|Zu \partial v|+|\partial u Z v| \right).
          \end{split}  \label{add 2.2}
       \end{equation} 
 where $Q_{ab}(u,v)=\partial_a u\partial_bv-\partial_av\partial_bu$, $a,b\in\{0,1,2,3\}$.
\end{lemma}
For the good derivatives, we have the following estimate
\begin{lemma}\cite{hou2020global} For  sufficiently regular functions  $u$, we have 
    \begin{equation}
|G_{a}u| \lesssim \langle t + |x| \rangle^{-1} ( |S u| + | Z u| ).
\end{equation}
\end{lemma}
\subsection{Sobolev-type inequalities}
\begin{lemma}\cite{georgiev1992decay}For sufficiently nice function $f$, we have the following Sobolev-type inequalities,\\
       (\romannum{1})(Klainerman-Sobolev inequality) 
        \begin{equation}
            |f(t,x)|\lesssim \langle t+r \rangle^{-1}\langle t-r \rangle^{-\frac{1}{2}}\sum_{|I|\le 2}|| \Gamma ^{I} f||_{L^2}.
        \end{equation}
    (\romannum{2})(Standard Sobolev inequality) 
        \begin{equation}
            |f(t,x)|\lesssim \langle r \rangle^{-1}\sum_{|I|\le 2}|| Z ^{I} f||_{L^2}.
        \end{equation}
         (\romannum{3})(Estimate inside the light cone) For $|x|\le \frac{t}{2}$,we have 
        \begin{equation}
            |f(t,x)|\lesssim \langle t \rangle^{-\frac{3}{4}}\sum_{|I|\le 3}|| Z ^{I} f||_{L^2}.
            \end{equation}
            (\romannum{4})(Global Sobolev inequality) 
        \begin{equation}
            |f(t,x)|\lesssim \langle t +r\rangle^{-\frac{3}{4}}\sum_{|I|\le 3}|| Z ^{I} f||_{L^2}. \label{2.6}
            \end{equation}
\end{lemma}
\begin{lemma}
   \cite{lindblad1990lifespan}  (Hardy's inequality) For a sufficiently nice function $f$, we have 
    \begin{equation}
        \left\|\frac{f}{|x|} \right \|_{L^2} \lesssim \left\| \nabla f\right\|_{L^2}.
    \end{equation}
\end{lemma}
Using H\"{o}lder's inequality, we have the following corollary.
\begin{corollary}
    Suppose $f$ is sufficiently nice, then for $1>\delta>0$, we have 
    \begin{equation}
         \left\|\frac{f}{|x|^\delta} \right \|_{L^2} \lesssim  \|\nabla f\|^{\delta}_{L^2} \|f\|^{1-\delta}_{L^2}.
    \end{equation}
\end{corollary}
\begin{proof}
    Write $\|\frac{f}{|x|^\delta}\|^2_{L^2}=\|\frac{f^2}{|x|^{2\delta}}\|^{2}_{L^1}$, by  H\"{o}lder's inequality, and Hardy's inequality,
    \begin{equation*}
   \begin{split}
         &  \|(\frac{f^{2\delta}}{|x|^{2\delta}})  f^{2-2\delta}\|_{L^1} \le \|\frac{f^{2\delta}}{|x|^{2\delta}}\|_{L^{\frac{1}{\delta}}} 
        \| f^{2-2\delta}\|_{L^{\frac{1}{1-\delta}}}
        \\ &= \|\frac{f^{2}}{|x|^{2}}\|_{L^1} ^{\delta} \|{f^{2}}\|_{L^1} ^{1-\delta} \lesssim \|\nabla f\|_{L^2}^2 \|f\|_{L^2}^2.
   \end{split}
    \end{equation*}
    Square root both sides, we obtain the corollary.
\end{proof}
\begin{remark}
     To the best of the author's knowledge, this corollary was first used in \cite{dong2023cubic} to tackle the problem caused by the largeness of the initial data.
\end{remark}
\subsection{Estimate for 3D wave equation }
We first state a weighted $L^\infty-L^\infty$ estimate for the linear wave equation. It's established in \cite{katayama2012global}.
 We define 
\begin{equation}
    W_\rho(t, r) := 
\begin{cases} 
\langle t + r\rangle^\rho & \text{if } \rho< 0, \\
\left[ \log \left( 2 + \langle t + r) \langle t - r\rangle^{-1} \right) \right]^{-1} & \text{if } \rho = 0, \\
\langle t - r\rangle^\rho & \text{if } \rho > 0.
\end{cases}
\end{equation}
And 
\[
W_-(t, r) := \min \{ \langle r \rangle , \langle t - r\rangle \}.
\]
\begin{lemma}
($L^\infty-L^\infty$ estimate)Let $w$ be a smooth solution to
\[
\Box w(t, x) = \Psi(t, x), \quad (t, x) \in (0, \infty) \times \mathbb{R}^3
\]
with initial data $w = \partial_t w = 0$ at $t = 0$.

Suppose that $\rho \geq 0$, $\kappa \geq 1$, and $\mu > 0$. Then there exists a positive constant $C = C(\rho, \kappa, \mu)$ such that
\begin{align*}
   & \langle t + |x|\rangle^{1-\rho} W_{\kappa-1}(t, |x|) |w(t, x)|
  \\ & \leq    C \sup_{\substack{\tau \in [0, t] \\ |y| - |x| \leq    t - \tau}} |y| \langle \tau+|y|\rangle^{\kappa + \mu-\rho} W_-(\tau, |y|)^{1-\mu} |\Psi(\tau, y)| 
   \\ & \leq    C \sup_{\substack{\tau \in [0, t] \\ y\in \mathbb{R}^3}  } \langle \tau+|y|\rangle^{2+\kappa  -\rho}  |\Psi(\tau, y)| 
\end{align*}
for $(t, x) \in (0, \infty) \times \mathbb{R}^3$.
\end{lemma}
Next, we state the standard conformal energy estimate and the standard energy estimate.

Given a sufficiently nice function $w = w(t, x)$, we define its energy on the constant time slice $t = constant$ by
\begin{equation}
   \mathcal{E} (t, w)
:=
\int_{\mathbb{R}^3} \Big(|\partial_t w|^2 + \sum_a |\partial_a w|^2 \Big) \, dx.
\end{equation}
And its conformal energy 
\begin{equation}
\begin{split}
      \mathcal{F} (t, w)
:&=
\int_{\mathbb{R}^3} \Big(| w|^2   +|\Omega w|^2+|L w|^2 +|Sw|^2\Big) \, dx 
\\ & =
\int_{\mathbb{R}^2} |  w|^2   +\sum_{a<b}|\Omega_{ab}  w|^2+\sum_a|L_a  w|^2 +|S w|^2 \, dx.
\end{split} 
\end{equation}
\begin{lemma}
    \cite{sogge1995lectures} Let $w$ be a smooth solution to
\[
\Box w(t, x) = \Psi(t, x), \quad (t, x) \in (0, \infty) \times \mathbb{R}^3
\]
with initial data $ w = w_0$, $\partial_t w = w_1$ at $t = 0$.Then we have the following\\
(\romannum{1}) (Conformal energy estimate) 
\begin{equation}
    \mathcal{F}^{\frac{1}{2}}(t, w) \lesssim \mathcal{F}^{\frac{1}{2}}(0, w) + \int_{0}^{t} \left\| \langle \tau + r \rangle \Psi(\tau, x) \right\|_{L^{2}} \, d\tau.
\end{equation}
(\romannum{2})(Standard energy estimate)
\begin{equation}
    \mathcal{E}^{\frac{1}{2}}(t, w) \lesssim \mathcal{E}^{\frac{1}{2}}(0, w) + \int_{0}^{t} \left\|  \Psi(\tau, x) \right\|_{L^{2}} \, d\tau.
\end{equation}

\end{lemma}
\subsection{Estimates for 3D Klein-Gordon equation}
We first state a $L^\infty-L^\infty$ estimate for the Klein-Gordon equation.
\begin{lemma}
   \cite{dong2024global} Let $u$ be the solution to the Cauchy problem
\[
-\Box u + u = G \quad \text{with} \quad (u, \partial_t u)_{t=0} = (u_0, u_1).
\]
Then we have
\begin{equation}
\begin{split}
   & \langle t + r\rangle^{\frac{3}{2}} |u(t, x)| \\ &  \lesssim  \sum_{k=0}^\infty \sum_{|l| \leq 5} \left\| |x|^{\frac{3}{2}} \zeta_k(|x|) Z^l u(0, x) \right\|_{L^2}
 \\ &+ \sum_{k=0}^\infty \sum_{|l| \leq 4} \max_{0 \leq s \leq t} \zeta_k(\tau) \left\| (\tau + |x|) Z^l G(\tau, x) \right\|_{L^2}.
\end{split}
\end{equation}
\end{lemma}
From the lemma above and the Cauchy-Schwarz inequality, we have 
 a simplified version of the lemma.
\begin{corollary}
  \cite{dong2024global}  In the context of Lemma 2.8, let $\delta_0 > 0$ and assume
\[
\sum_{|I| \leq 4} \max_{2 \leq \tau \leq t} \tau^{\delta_0} \left\| (\tau + |x|) Z^I G(\tau, x) \right\|_{L^2} \leq C_G.
\]

Then we have
\[
\langle t + r\rangle^{\frac{3}{2}} |u(t, x)| \lesssim \frac{C_G}{1 - 2^{- \delta_0}} + \sum_{|I| \leq 5} \|  \langle|x|\rangle ^{\frac{3}{2}} \log(2 + |x|) Z^I u(0, x) \|_{L^2}.
\]
\end{corollary}
Then we state the standard energy estimate and the ghost weight estimate for the Klein-Gordon equation. Their proof can be found in \cite{ alinhac2009hyperbolic, dong2021global}.
\begin{lemma}(Standard energy estimate)
   Assume $w$ is the solution to
\[
-\Box w +  w = f.
\]
Then we have
\begin{align}
E^{\frac{1}{2}}(t,w) &\lesssim E^{\frac{1}{2}}(0,w) +  \int_0^t \|f\|_{L^2} (\tau)d\tau, \tag{2.10}
\end{align}
where  \begin{equation}
 {E} (t, w)
:=
\int_{\mathbb{R}^3} \Big(|\partial_t w|^2 + \sum_a |\partial_a w|^2+|w|^2 \Big) \, dx.
\end{equation}
\end{lemma}
\begin{lemma}
    (Ghost weight estimate) Assume $w$ is the solution to
\begin{equation}
    -\Box w +  w = f.
\end{equation}
Then we have
\begin{equation}
    \begin{split}
       E_{gst}(t,w) &\lesssim\int_{\mathbb{R}^3}  (|\partial_t w|^2 + \sum_a |\partial_a w|^2 +  w^2) \, dx(0)\\& +  \int_0^t \int_{\mathbb{R}^3} |f \partial_t w | \, dxdt,
    \end{split}
\end{equation}
\begin{comment}
    \begin{equation*}
    q = \int_{-\infty}^{r-t} \langle s \rangle ^{-3/2} \, ds,
\end{equation*}
and
\end{comment}
in which
\begin{equation}
    \begin{split}
        E_{gst}(t,w) &= \int_{\mathbb{R}^3}  (|\partial_t w|^2 + \sum_a |\partial_a w|^2 +  w^2) \, dx(t) \\ & \int_0^t \int_{\mathbb{R}^3} \frac{w^2}{\langle r-t \rangle ^{9/8}}  \, dxdt
+ \sum_a \int_0^t \int_{\mathbb{R}^2} \frac{|G_{a}w|^2}{\langle r-t \rangle^{9/8}}  \, dxdt.
    \end{split}
\end{equation}
\end{lemma}
\section{Proof of Theorem 1} \label{sec3}
\subsection{Bootstrap assumptions} \label{sec3.1}
Let $N>40$ be an integer, $0<p<\frac{1}{10}$, we consider the following bootstrap assumptions: For $C\ge10$ and $0<\epsilon<C^{-1}$ to be chosen later,
we have
\begin{equation}
\left\{
\begin{aligned}
&\sum_{|I| \leq N}   \mathcal{E}^{\frac{1}{2}} (t,Z^Iu)\leq C \langle t\rangle^p, \\
&\sum_{|I| \leq N-1}    \mathcal{E}^{\frac{1}{2}} (t,\partial Z^Iu)\leq C \epsilon\langle t\rangle^p, \\
&\sum_{|I| \leq N-10} \mathcal{F}^{\frac{1}{2}}(t,Z^Iu)  \leq \frac{C}{\epsilon}\langle t\rangle^p ,
\\
&\sum_{|I| \leq N-10} \mathcal{F}^{\frac{1}{2}}(t,\partial Z^Iu)  \leq C \langle t\rangle^p,
\\
&\sum_{|I| \leq N} {E}^{\frac{1}{2}}_{gst}(t, Z^Iw)  \leq C \epsilon  .
 \label{3.1}
\end{aligned}
\right.
\end{equation}
For initial data $(u_0,u_1,w_0,w_1)$, we set
\begin{equation}
    T_{*}:= T_{*}(u_0,u_1,w_0,w_1)=\sup \left\{ {t\in[0,\infty)}\big|(u,w) \text{  satisfying \eqref{3.1}on [0,$t$)}\right\}.
\end{equation}
We have the following proposition.
\begin{proposition}
    For any initial data $(u_0,u_1,w_0,w_1)$ satisfying the requirements \eqref{1.3} in Theorem 1.1, we have $T_{*}=\infty$.
\end{proposition}
The initial value automatically satisfies \eqref{3.1} under the requirements of \eqref{1.3}; the proof is based on the equations and an elementary induction argument on the size of the index. We list all the needed results without a detailed proof.
\subsection{Estimate of initial values}
\begin{claim}
     For $|I|\le N-8$, we have 
    \begin{equation}
          \sum_{|J| \leq 5} \|  \langle |x| \rangle ^{\frac{3}{2}} \log(2 + |x|) Z^J Z^I w(0,x) \|_{L^2} \le 10 \epsilon.
    \end{equation}
\end{claim}
\begin{comment}
     \begin{proof}
       Using $\langle|x|\rangle^{\frac{3}{2}}|log (2+|x|)|\le \langle|x|\rangle^{\frac{5}{2}}+\langle|x|\rangle^{\frac{3}{2}}\le2 \langle x\rangle^{\frac{5}{2}}$ , we have 
       \begin{equation*}
       \begin{split}
           & \sum_{|J| \leq 5} \| \langle |x| \rangle^{\frac{3}{2}} \log(2 + |x|) Z^J Z^I w(0,x) \|_{L^2}
           \\ &\le    2 \sum_{|I| \leq N-3} \|  \langle x\rangle ^{\frac{5}{2}+|I|}  \nabla^I w_0 \|_{L^2}+
           \\ &\le 2\epsilon.
       \end{split} 
       \end{equation*}
   \end{proof}
\end{comment}
\begin{claim}
  For $|I|\le N-1$, we have 
  \begin{equation}
  \begin{split}
  \mathcal{E}^{\frac{1}{2}}(0,\partial Z^Iu) \le 2\epsilon.
  \end{split}
  \end{equation}
  In particular, for $|I|\le N-2$, we have 
  \begin{equation}
  \begin{split}
  \mathcal{E}^{\frac{1}{2}}(0,\partial \partial Z^Iu) \le 2\epsilon.
  \end{split}
  \end{equation}
  And also, for $|I|\le N$,
  \begin{equation}
 \mathcal{E}^{\frac{1}{2}}(0, Z^Iu) \le 2.   
  \end{equation}
\end{claim}
\begin{comment}
    \begin{proof}
   Note that $\|\nabla Z^I u(0,\cdot)\|_{\dot {H}^1}=\|\nabla \nabla  Z^I u(0,\cdot)\|_{L^2}\le\|\langle|x|\rangle^{|I|}\nabla  \nabla  \nabla ^I u_0\|_{L^2}$, from \ref{1.3} we have $\|\nabla   Z^I u_0\|_{\dot {H}^1}\le \epsilon$. Following the same manner, we have $\| \nabla  Z^I u_1\|_{L^2}\le \epsilon$, therefore we derive the first inequality, the others follow from the same manner.
\end{proof}
\end{comment}
\begin{claim}
For $|I|\le N-10$, we have 
\begin{equation}
\begin{split}
     \mathcal{F}^{\frac{1}{2}}(0,   Z^I u) \le\frac{200}{\epsilon},\quad
       \mathcal{F}^{\frac{1}{2}}(  0, \partial Z^I u) \le 10,
       \quad \mathcal{F}^{\frac{1}{2}}(  0, S Z^I u)\le 10.
\end{split} 
\end{equation}
\end{claim}
\begin{comment}
    \begin{proof}
   From the definition of conformal energy, we have 
   \begin{equation}
   \begin{split}
     \mathcal{F}^{\frac{1}{2}}(0,   Z^I u)&=
\Big(\int_{\mathbb{R}^2} | Z^I u|^2   +\sum_{a<b}|\Omega_{ab} Z^I u|^2+\sum_a|L_a Z^I u|^2 +|SZ^I u|^2 \, dx\Big)^{\frac{1}{2}}  
\\ & \le \|Z^I u_0\|_{L^2}+\sum_{a<b}\|\Omega_{ab} Z^I u_0\|_{L^2}+\sum_a\|L_a Z^I u_0\|_{L^2}+\|SZ^I u_0\|_{L^2}
\\ &\le  \|\langle |x|\rangle ^{|I|}\partial^I u_0\|_{L^2}+ 7\|\langle |x|\rangle ^{|I|+1}\partial\partial^I u_0\|_{L^2}
\\ & \le \frac{100}{\epsilon}+ 7 \le \frac{200}{\epsilon}.
   \end{split}
       \end{equation}
       The other inequalities follow from the same approach.
\end{proof}
\end{comment}
\begin{claim}
For $|I|\le N$, we have 
\begin{equation}
    E_{gst}^{\frac{1}{2}}(0,Z^Iw) \le{10} \epsilon.
\end{equation}
\end{claim}
\begin{comment}
    \begin{proof}
    It follows directly from \ref{1.3}.
\end{proof}
\end{comment}

We then derive estimates for lower index terms under bootstrap assumptions.
\subsection{Estimates under bootstrap assumptions}
Using Klainerman's inequality directly, we immediately obtain the following.
\begin{proposition}
    Let $|I|\le N-3$, for all $\tau\in[0,T_*)$, we have 
    \begin{equation}
       \begin{split}
         &  | \partial Z^I  u(\tau,x)|\lesssim C  \langle \tau +r\rangle^{-\frac{3}{4}+p},
         \\ &  | Z^I  w(\tau,x)|\lesssim C \epsilon  \langle \tau +r\rangle^{-\frac{3}{4}}.
       \end{split} 
    \end{equation} 
    For $|I|\le N-13$, for all $\tau\in[0,T_*)$, we have
    \begin{equation}
       |SZ^I  u (\tau,x)|+|ZZ^I  u (\tau,x)|\lesssim  \frac{C}{\epsilon} \langle \tau +r\rangle^{-\frac{3}{4}+p}.
    \end{equation} 
\end{proposition}
\begin{proof}
    Using \eqref{lemma 2.1},  \eqref{2.6} and \eqref{3.1}, we have 
\begin{equation*}
  \sum_{|I|\le N-3} |Z^I \partial u(\tau,x)|\lesssim  \sum_{|I|\le N-3} \sum_{|J|\le 3} ||Z^I Z^J \partial u||_{L^2} \langle \tau +r \rangle ^{-\frac{3}{4}}\le C\epsilon  \langle \tau +r \rangle ^{-\frac{3}{4}+p}.
\end{equation*}
The other inequalities follow from the same approach.
\end{proof}
Next, we derive a point-wise estimate for the Klein-Gordon equation.
\begin{proposition}
     For $|I|\le N-14$,  $\tau\in[0,T_*)$, we have
     \begin{equation}
         |Z^I w(\tau,x)| \lesssim C^2 \epsilon \langle \tau+r \rangle ^{-\frac{3}{2}}.
     \end{equation}
\end{proposition}
\begin{proof}
    For $|I|\le N-14$, act $ Z^I$ on Klein-Gordon equation in \eqref{add 1.2}, using Lemma 2.1, we have
    \begin{equation}
        -\Box  Z^Iw+ Z^Iw= Z^IQ_0(u,w).
    \end{equation}
In view of Corollary 2.10, we have to control
\begin{equation}
   \sum_{|I_1|+|I_2|\le N-10} ||\langle \tau +r\rangle Q_0(Z^{I_1}u,Z^{I_2}w)||_{L^2}\label{3.6}.
\end{equation}
Using Lemma 2.2 and H\"{o}lder's inequality, we have 
\begin{equation}
\begin{split}
     &  \sum_{|I_1|+|I_2|\le N-10} ||\langle \tau +r\rangle Q_0(Z^{I_1}u,Z^{I_2}w)||_{L^2}   
     \\  &\lesssim
        \sum_{|I_1|+|I_2|\le N-10} \|(|S Z^{I_1} u|+|ZZ^{I_1} u|) Z Z^{I_2} w\|_{L^2} 
        \\  &\lesssim
        \sum_{|I_1|+|I_2|\le N-10} (||S Z^{I_1} u||_{L^6}+||ZZ^{I_1} u||_{L^6} )||ZZ^{I_2} w||_{L^3}.
\end{split}
    \label{3.8}
\end{equation}
For $L^6$ part above, using Sobolev embedding $\iota: \dot H^1(\mathbb{R}^3)\rightarrow L^6(\mathbb{R}^3)$, we have 
\begin{equation}
\begin{split}
 &  \sum_{|I_1|\le N-10} ||S Z^{I_1} u||_{L^6}  +||ZZ^{I_1} u||_{L^6}   
  \\ &   \lesssim   \sum_{|I_1|\le N-10} ||\partial S  Z^{I_1} u||_{L^2} +||\partial ZZ^{I_1} u||_{L^2}   
   \\ &  \lesssim  \sum_{|I_1|\le N-10} \mathcal{F}^{\frac{1}{2}}(\tau,\partial Z^{I_1}u)+  \sum_{|I_1|\le N-1} ||\partial Z^{I_1} u||_{L^2}
   \\ & \lesssim C \langle \tau\rangle ^p \label{3.9}.
\end{split}
\end{equation}
For $L^3$ part in \eqref{3.8}, using $L^p$ interpolation, by \eqref{3.1} and Proposition 3.2, we have 
\begin{equation}
\left\{
\begin{aligned}
& \sum_{|I_2|\le N-10} || Z  Z^{I_2} w||_{L^\infty} \lesssim C \epsilon \langle \tau \rangle ^{-\frac{3}{4}},  \\ &
 \sum_{|I_2|\le N-10} || Z  Z^{I_2} w||_{L^2} \lesssim C \epsilon. \
   \end{aligned}
   \right.
\end{equation}
Therefore, interpolating $L^3$ by $L^2$ and $L^\infty$, we have \begin{equation}
    \sum_{|I_2|\le N-10} || Z  Z^{I_2} w||_{L^3} \lesssim C \epsilon \langle \tau \rangle ^{-\frac{1}{4}}. \label{3.10}
\end{equation}
Substituting \eqref{3.8} and \eqref{3.10} into \eqref{3.6}, since $0<p<\frac{1}{10}$, we have
\begin{equation}
\begin{split}
      \sum_{|I_1|+|I_2|\le N-10} ||\langle \tau +r\rangle Q_0(Z^{I_1}u,Z^{I_2}w)||_{L^2} \lesssim C^2 \epsilon \langle \tau \rangle ^{-\frac{1}{4}+p}\le C^2 \epsilon \langle \tau \rangle ^{-\frac{1}{10}}.
\end{split}
\end{equation}

Therefore, we may choose $\delta_0=\frac{1}{10}$ in Corollary 2.10, combining Claim 1, we have 
\begin{equation}
   \langle \tau+r\rangle^{\frac{3}{2}} |Z^I w(\tau,x)| \lesssim C^2 \epsilon +\epsilon \lesssim C^2 \epsilon.
\end{equation}
That's the desired result since we require $C>10$.
\end{proof}
    \begin{proposition}
     Let $|I|\le N-17$, for all $\tau\in[0,T_*)$, $0\le\lambda\le \frac{1}{100}$, we have 
    \begin{equation}
        | \partial Z^Iu(\tau,x)|\lesssim_{\lambda} C^{3}\epsilon^{\lambda ^2} \langle \tau +r\rangle^{-1+\frac{3\lambda}{2}} \langle \tau -r\rangle^{-\frac{1}{8}}.
    \end{equation} 
\end{proposition}
\begin{proof}
    Acting $\partial\partial Z^I$ over the wave equation in \eqref{add 1.2}, we have
    \begin{equation}
         -\Box  \partial \partial Z^Iu= \partial \partial Z^IQ_0(u,w).
    \end{equation}

We write $\partial\partial Z^I u=\hat{u}+\bar{u}$, where $\hat{u}$ solves
\begin{equation}
\left\{
\begin{aligned}
& -\Box \hat{u}=\partial \partial Z^IQ_0(u,w) \\ &
 (\hat{u}(0),\partial_t\hat{u}(0))=(0,0)
  \end{aligned}
   \right.
\end{equation}
We first derive a point-wise estimate for $\hat{u}$. In view of Lemma 2.6, it suffices to pointwisely control $|\partial \partial Z^IQ_0(u,w)|$. Applying Lemma 2.1, Lemma 2.2, Proposition 3.2, and Proposition 3.3, we have the following
\begin{equation}
    \begin{split}
        |\partial \partial Z^I Q_0(u,w)| & \lesssim \sum_{|I_1|+|I_2|\le N-15} |Q_0(Z^{I_1}u,Z^{I_2}w)|
        \\ & \lesssim \sum_{|I_1|+|I_2|\le N-15}\frac{1}{\langle \tau+r\rangle} (|S Z^{I_1}u|+|ZZ^{I_1}u|)|ZZ^{I_2}w|
       \\ &   \lesssim \frac{1}{\langle \tau+r\rangle} \frac{C}{\epsilon} {\langle \tau+r\rangle}^{-\frac{3}{4}+p} C^2 \epsilon  {\langle \tau+r\rangle}^{-\frac{3}{2}}
       \\ & = C^3  {\langle \tau+r\rangle}^{-\frac{13}{4}+p}.
    \end{split}
\end{equation}
Choose $\rho=0$,$\kappa=\frac{5}{4}-p$ in Lemma 2.7, we have 
\begin{equation}
    |\hat{u}(t,x)|\lesssim C^3 {\langle \tau+r\rangle}^{-1}{\langle \tau-r\rangle}^{-\frac{1}{4}+p} \label{3.17}.
\end{equation}
Next, we derive a point-wise estimate for $\bar{u}$, where $\bar{u}$ solves liner wave equation
\begin{equation}
\left\{
\begin{aligned}
& -\Box \bar{u}=0 \\ &
 (\bar{u}(0),\partial_t\bar{u}(0))=(\partial \partial Z^I u(0,\cdot),\partial_t\partial \partial Z^I u(0,\cdot))
   \end{aligned}
   \right.
\end{equation}
Based on Claim 2 and Klainerman's inequality, we have 
\begin{equation}
      |\bar{u}(t,x)|\lesssim \epsilon{\langle \tau+r\rangle}^{-1}{\langle \tau-r\rangle}^{-\frac{1}{2}} \label{3.19}.
\end{equation}
Using \eqref{3.17} and \eqref{3.19}, we have point-wise estimate for $\partial \partial Z^I u$
\begin{equation}
     |\partial \partial Z^I {u}(t,x)|\le   |\bar{u}(t,x)|+ |\hat{u}(t,x)|\lesssim C^3 {\langle \tau+r\rangle}^{-1}{\langle \tau-r\rangle}^{-\frac{1}{4}+p} \label{3.20}.
\end{equation}
After repeating the almost same process, for the same requirements of index $I$, we can also obtain
\begin{equation}
     |\partial  Z^I {u}(t,x)|\le   |\bar{u}(t,x)|+ |\hat{u}(t,x)|\lesssim C^3 {\langle \tau+r\rangle}^{-1}{\langle \tau-r\rangle}^{-\frac{1}{4}+p}. \label{3.21}
\end{equation}
By Proposition 3.2, we also have 
\begin{equation}
    | \partial \partial Z^I  u(\tau,x)|\lesssim C  \epsilon \langle \tau +r\rangle^{-\frac{3}{4}+p} \label{3.22}.
\end{equation}
Interpolating \eqref{3.20} and \eqref{3.22}, using $0<p<\frac{1}{10}$, for $\lambda\in[0,1]$, we have 
\begin{equation}
\begin{split}
       | \partial \partial Z^I  u(\tau,x)|& \lesssim C^{3-2\lambda}  \epsilon ^{\lambda} \langle \tau +r\rangle^{-1+\lambda p+\frac{\lambda}{4}} \langle \tau -r\rangle^{-\frac{1}{4}+ p+\frac{\lambda}{4}-p\lambda}
       \\ & \le C^{3-2\lambda}  \epsilon ^{\lambda} \langle \tau +r\rangle^{-1+ \frac{\lambda}{2}} \langle \tau -r\rangle^{-\frac{1}{8}+ \frac{\lambda}{4}}.
\end{split}
\end{equation}
Integrating along the radial direction, since $\partial Z^I u$ vanishes at infinity, we have 
\begin{equation}
\begin{split}
     |\partial Z^Iu(t,x)| & \le \int_{|x|}^{\infty}|\partial_r \partial Z^I u(t,y)|dy 
    \\ & \lesssim  \int_{|x|}^{\infty} C^{3-2\lambda}  \epsilon ^{\lambda} \langle \tau +r\rangle^{-1+ \frac{\lambda}{2}} \langle \tau -r\rangle^{-\frac{1}{8}+ \frac{\lambda}{4}} dr
    \\ & \le  C^{3-2\lambda}  \epsilon ^{\lambda} \int_{|x|}^{\infty}  {\langle \tau+r\rangle}^{-\frac{1}{8}+\lambda}{\langle \tau-r\rangle}^{-1-\frac{\lambda}{4}} dr
    \\ & \lesssim  C^{3-2\lambda}  \epsilon ^{\lambda}  {\langle \tau+r\rangle}^{-\frac{1}{8}+\lambda} \label{3.24}.
\end{split}
\end{equation}
Interpolating \eqref{3.21} and \eqref{3.24}, and also using $0\le\lambda\le \frac{1}{100}$, we finally obtain 
\begin{equation}
\begin{split}
       |\partial Z^Iu(t,x)| &\lesssim  C^{3-2\lambda^2}\epsilon^{\lambda^2} {\langle \tau+r\rangle}^{-1+\frac{7\lambda}{8}+\lambda^2}{\langle \tau-r\rangle}^{-\frac{1}{4}+\frac{1}{\lambda}+p(1-\lambda)}
       \\ &\lesssim C^{3}\epsilon^{\lambda^2} {\langle \tau+r\rangle}^{-1+\frac{3\lambda}{2}}{\langle \tau-r\rangle}^{-\frac{1}{8}}.
\end{split}
\end{equation}
\end{proof}
Next, we state the last necessary point-wise estimate for $w$.
\begin{proposition}
    For $|I|\le N-8$, $\tau\in[0,T_*)$, we have 
    \begin{equation}
        |Z^I w(\tau ,x)| \lesssim  C^4 \epsilon^{\frac{5}{7}}\langle \tau +r \rangle^{-\frac{3}{2}}.
    \end{equation}
\end{proposition}
\begin{proof}
    Like the proof of Proposition 3.3, in view of Corollary 2.10, it suffices to control 
    \begin{equation}
        \sum_{|I_1|+|I_2|\le N-4} ||\langle \tau+r \rangle Q_0(Z^{I_1}u,Z^{I_2}w)||_{L^2}.
    \end{equation}
    We split the term above as follows
     \begin{equation}
      \begin{split}
         &  \sum_{|I_1|+|I_2|\le N-4} ||\langle \tau+r \rangle Q_0(Z^{I_1}u,Z^{I_2}w)||_{L^2} \\ &\lesssim  \sum_{|I_1|+|I_2|\le N-4}   ||(|S Z^{I_1} u|+|ZZ^{I_1} u|) Z Z^{I_2} w||_{L^2} 
         \\ & \lesssim \sum_{\substack{|I_1|\le N-4 \\ |I_2|\le N-15}}   ||(|S Z^{I_1} u|+|ZZ^{I_1} u|) Z Z^{I_2} w||_{L^2} 
         \\ & + \sum_{\substack{|I_1|\le N-17 \\ |I_2|\le N-4}} ||(|S Z^{I_1} u|+|ZZ^{I_1} u|) Z Z^{I_2} w||_{L^2} .
      \end{split} 
    \end{equation}
    For the first term, using \eqref{3.1}, Lemma 2.2 and Proposition 3.3, we have 
    \begin{equation}
        \begin{split}
            &  \sum_{\substack{|I_1|\le N-4 \\ |I_2|\le N-15}}   ||(S Z^{I_1} u+ZZ^{I_1} u) Z Z^{I_2} w||_{L^2} 
            \\ & \lesssim \sum_{\substack{|I_1|\le N-4 \\|I_2|\le N-15}}   ||\langle \tau +r\rangle \partial Z^{I_1} u Z Z^{I_2} w||_{L^2} 
            \\ & \lesssim \sum_{\substack{|I_1|\le N-4 \\ |I_2|\le N-15}}  || \partial Z^{I_1} u ||_{L^2} ||\langle \tau +r \rangle Z Z^{I_2} w ||_{L^\infty} 
           \\ & \lesssim C^3 \epsilon \langle\tau\rangle ^{-\frac{1}{2}+p}. \label{3.29}
        \end{split}
    \end{equation}
    For the second term, we have 
    \begin{equation}
    \begin{split}
       &  \sum_{\substack{|I_1|\le N-17 \\|I_2|\le N-4}}   \|(|S Z^{I_1} u|+|ZZ^{I_1} u|) Z Z^{I_2} w\|_{L^2}  
       \\ & \le \sum_{\substack{|I_1|\le N-17 \\ |I_2|\le N-4}}   \||S Z^{I_1} u|+|ZZ^{I_1} u|\|_{L^7}  || Z Z^{I_2} w||_{L^\frac{14}{5}} \label{3.30}.
    \end{split}   
    \end{equation}
    For the $L^7$ term above, using \ref{3.1} and Proposition 3.4, we have the following 
    \begin{equation}
    \begin{split}
      & ||S Z^{I_1} u||_{L^2}+||ZZ^{I_1} u||_{L^2}\lesssim \frac{C}{\epsilon}\langle \tau \rangle^{p},
      \\ & ||S Z^{I_1} u||_{L^\infty}+||ZZ^{I_1} u||_{L^\infty} \\ &\lesssim ||\langle \tau +r \rangle\partial Z^{I_1} u||_{L^\infty}
      \\ & \lesssim C^{3}\epsilon^{\lambda^2} {\langle \tau\rangle}^{\frac{3\lambda}{2}}.
      \end{split}
\end{equation}
Interpolating $L^7$ by $L^2$ and $L^\infty$, we derive that
\begin{equation}
      ||S Z^{I_1} u||_{L^7}+||ZZ^{I_1} u||_{L^7}\lesssim C^{\frac{17}{7}}\epsilon^{-\frac{2}{7}+\frac{5\lambda^2}{7}} \langle \tau  \rangle^{\frac{15\lambda}{14}+\frac{2p}{7}}. \label{3.32}
\end{equation}
For the $L^{\frac{14}{5}}$ part, from \eqref{3.1} and Proposition 3.2 we have 
\begin{equation}
\begin{split}
    &  ||ZZ^I  w (\tau,x)||_{L^{\infty
    }}\lesssim C \epsilon \langle \tau \rangle^{-\frac{3}{4}},
    \\ & ||ZZ^I  w (\tau,x)||_{L^{2
    }}\lesssim C \epsilon .
\end{split}
\end{equation}
Interpolating $L^{\frac{14}{5}}$ by $L^2$ and $L^\infty$, we derive that
\begin{equation}
  ||ZZ^I  w (\tau,x)||_{L^{\frac{14}{5}
    }}\lesssim C \epsilon \langle \tau \rangle^{-\frac{3}{14}} .\label{3.34}
\end{equation}
Using \eqref{3.32} and \eqref{3.34}, substitute into \eqref{3.30}, we obtain that
\begin{equation}
    \begin{split}
       &  \sum_{\substack{|I_1|\le N-17 \\ |I_2|\le N-4}}  ||(|S Z^{I_1} u|+|ZZ^{I_1} u|) Z Z^{I_2} w||_{L^2}  
       \\ & \lesssim  
       C^{\frac{24}{7}}\epsilon^{\frac{5}{7}+\frac{5\lambda^2}{7}} \langle \tau  \rangle^{\frac{15\lambda}{14}-\frac{3}{14}+\frac{2p}{7}} 
       \\ & \lesssim C^4 \epsilon ^{\frac{5}{7}} \langle \tau \rangle^{-\frac{1}{14}} \label{3.35}.
    \end{split}   
    \end{equation}
Using \eqref{3.29} and \eqref{3.35} , we finally derive that 
\begin{equation}
\begin{split}
       \sum_{|I_1|+|I_2|\le N-4} ||\langle \tau+r \rangle Q_0(Z^{I_1}u,Z^{I_2}w)||_{L^2} & \lesssim  C^3 \epsilon \langle\tau\rangle ^{-\frac{1}{4}}+ C^4 \epsilon ^{\frac{5}{7}} \langle \tau \rangle^{-\frac{1}{14}}
       \\ & \lesssim  C^4 \epsilon ^{\frac{5}{7}} \langle \tau \rangle^{-\frac{1}{14}}
\end{split}
    \end{equation}
    Using  Corollary 2.10, choose $\delta_0=\frac{1}{14}$, and Claim 1 again, we obtain the desired result.
  \end{proof}  
    \begin{proposition}
        For $|I|\le N-16$, $t\in[0,T_*)$, we have
        \begin{equation}
            \mathcal{F}^{\frac{1}{2}}(t,SZ^Iu)\lesssim   C^3
            \epsilon^{\frac{1}{8}}\langle\tau\rangle^{\frac{5}{8}+p}+10.
        \end{equation}
    \end{proposition}
    \begin{proof}
       First we act $SZ^I$ on the wave equation, using \eqref{2.1},  
       \begin{equation}
       \begin{split}
                &   -\Box SZ^I u
                \\&=-2\Box Z^Iu+S(-\Box Z^I u)
                \\&=-2\Box Z^Iu+\sum_{|I_1|+|I_2|\le |I|}Q_0(SZ^{I_1}u,Z^{I_2}w)+2Q_0(Z^{I_1}u,SZ^{I_2}w)   
                   \\& =\sum_{|I_1|+|I_2|\le |I|} -2Q_0(Z^{I_1}u,Z^{I_2}w)+Q_0(SZ^{I_1}u,Z^{I_2}w)+Q_0(Z^{I_1}u,SZ^{I_2}w)
                   .
       \end{split}
       \end{equation}
To apply the conformal energy estimate(Lemma 2.8), we have the following :
\\ First, using Proposition 3.3, \ref{3.1}, we have 
\begin{equation}
\begin{split}
    & \sum_{|I_1|+|I_2|\le |I|} \|-2\langle\tau+r\rangle Q_0(Z^{I_1}u,Z^{I_2}w)\|_{L^2}(\tau)
    \\ & \lesssim \sum_{|I_1|+|I_2|\le N-16} \|\langle\tau+r\rangle \partial Z^{I_1}u\partial Z^{I_2}w\|_{L^2}(\tau)
    \\ & \lesssim \sum_{|I_1|+|I_2|\le N-16} C^2 \epsilon \langle \tau\rangle^{-\frac{1}{2}} \|\partial Z^{I_1}u\|_{L^2}(\tau)
    \\ & \lesssim C^3 \epsilon \langle \tau\rangle^{-\frac{1}{2}+p} \label{add 3.44}.
\end{split}
\end{equation}
Second, using Proposition 3.3, \eqref{3.1} we have 
\begin{equation}
    \begin{split}
         & \sum_{|I_1|+|I_2|\le |I|} \|\langle\tau+r\rangle Q_0(SZ^{I_1}u,Z^{I_2}w)\|_{L^2}(\tau)
    \\ & \lesssim \sum_{|I_1|+|I_2|\le N-16} \|\langle\tau+r\rangle \partial SZ^{I_1}u\partial   Z^{I_2}w\|_{L^2}(\tau)
    \\ & \lesssim \sum_{|I_1|\le N-16} C^2 \epsilon \langle \tau\rangle^{-\frac{1}{2}} \|\partial SZ^{I_1}u\|_{L^2}(\tau)
    \\ & \lesssim  C^3 \epsilon \langle \tau\rangle^{-\frac{1}{2}+p}\label{add 3.45}.
    \end{split}
\end{equation}
Third, using \ref{add 2.2}, \ref{3.1}, and Lemma 2.1, we have
\begin{equation}
     \begin{split}
         & \sum_{|I_1|+|I_2|\le |I|} \|\langle\tau+r\rangle Q_0(Z^{I_1}u,SZ^{I_2}w)\|_{L^2}(\tau)
    \\ & \lesssim \sum_{|I_1|+|I_2|\le N-16} \|\langle\tau+r\rangle \langle\tau+r\rangle^{-1}(|SZ^{I_1}u|+|ZZ^{I_1}u|)SZZ^{I_2}w\|_{L^2}(\tau)
    \\& + \|\langle\tau+r\rangle \langle\tau+r\rangle^{-1}(|SZ^{I_1}u|+|ZZ^{I_1}u|)\partial Z^{I_2}w\|_{L^2}(\tau)
    \\ & \lesssim \sum_{|I_1|+|I_2|\le N-16} \|\frac{SZ^{I_1}u}{r^{\frac{1}{8}}}\|_{L^2}\|r^{\frac{1}{8}}SZZ^{I_2}w\|_{L^{\infty}}+\|\frac{ZZ^{I_1}u}{r^{\frac{1}{8}}}\|_{L^2}\|r^{\frac{1}{8}}SZZ^{I_2}w\|_{L^{\infty}}
    \\ &+\|\frac{S Z^{I_1}u}{r^{\frac{1}{8}}}\|_{L^2}\|r^{\frac{1}{8}}\partial Z^{I_2}w\|_{L^{\infty}}+\|\frac{ZZ^{I_1}u}{r^{\frac{1}{8}}}\|_{L^2}\|r^{\frac{1}{8}}\partial Z^{I_2}w\|_{L^{\infty}} .
   \label{add 3.46}
    \end{split}
\end{equation}

It should be mentioned that the terms where $\partial$ hits $w$ are much better than the terms where $S$ hits $w$, because the vector field $S$ brings an extra $\langle \tau+r\rangle$ increase. Therefore, it only needs to treat the terms where $w$ is hit by $S$.

Using Corollary 2.6, \eqref{3.1},
\begin{equation}
    \begin{split}
     &   \|\frac{SZ^{I_1}u}{r^{\frac{1}{8}}}\|_{L^2}(\tau)\lesssim  \|\partial SZ^{I_1}u\|_{L^2}^{\frac{1}{8}} \| SZ^{I_1}u\|_{L^2}^{\frac{7}{8}} \lesssim C\epsilon^{-\frac{7}{8}}\langle \tau\rangle^{p},
     \\ &\|\frac{ZZ^{I_1}u}{r^{\frac{1}{8}}}\|_{L^2}(\tau)\lesssim  \|\partial ZZ^{I_1}u\|_{L^2}^{\frac{1}{8}} \| ZZ^{I_1}u\|_{L^2}^{\frac{7}{8}} \lesssim C\epsilon^{-\frac{7}{8}}\langle \tau\rangle^{p}. \label{ADD 3.47}
    \end{split}
\end{equation}
Using Proposition 3.3, 
\begin{equation}
\begin{split}
        \|r^{\frac{1}{8}}SZZ^{I_2}w\|_{L^{\infty}}\le \|r^{\frac{1}{8}}\langle \tau+r\rangle \partial ZZ^{I_2}w\|_{L^{\infty}}\lesssim C^2 \epsilon \langle \tau \rangle ^{-\frac{3}{8}}. \label{add 3.48}
\end{split}
\end{equation}
Substituting \eqref{ADD 3.47}, \eqref{add 3.48} into \eqref{add 3.46}, we have
\begin{equation}
     \begin{split}
         & \sum_{|I_1|+|I_2|\le |I|} \|\langle\tau+r\rangle Q_0(Z^{I_1}u,SZ^{I_2}w)\|_{L^2}(\tau)
     \lesssim C^3 \epsilon^{\frac{1}{8}} \langle \tau \rangle ^{-\frac{3}{8}+p}.
   \label{add 3.49}
    \end{split}
\end{equation}
Finally, using \ref{add 3.44},\ref{add 3.45}, \ref{add 3.49}, Claim 3, and the conformal energy estimate Lemma 2.8, we have the desired,
\begin{equation}
    \begin{split}
        & \mathcal{F}^{\frac{1}{2}}(SZ^Iu,t)
        \\ &\lesssim \mathcal{F}^{\frac{1}{2}}(SZ^Iu,0)+\int_0^t \sum_{|I_1|+|I_2|\le |I|} \|\langle\tau +r\rangle Q_0(Z^{I_1}u,Z^{I_2}w)\|_{L^2}(\tau)
        \\ &+\|\langle\tau +r\rangle Q_0(SZ^{I_1}u,Z^{I_2}w)\|_{L^2}(\tau)
         +\|\langle\tau +r\rangle Q_0(Z^{I_1}u,SZ^{I_2}w)\|_{L^2}(\tau)d \tau
        \\ &\lesssim \int_0^t C^3 \epsilon^{\frac{1}{8}} \langle \tau \rangle ^{-\frac{3}{8}+p}+C^3 \epsilon \langle \tau\rangle^{-\frac{1}{2}+p}+10
      \\&  \lesssim C^3 \epsilon^{\frac{1}{8}} \langle \tau \rangle ^{\frac{5}{8}+p}+10.
    \end{split}
\end{equation}
   \end{proof}
Then, we derive an estimate for the good derivatives of the wave components.
\begin{proposition}
    For $|I|\le N-19$, $0\le\lambda\le\frac{1}{100}$, $t\in[0,T_*)$, we have
    \begin{equation}
        |G_aZ^Iu(t,x)|\lesssim C^3 \epsilon^{\frac{\lambda^2}{2}}\langle t+r\rangle^{-\frac{17}{16}+\frac{p}{2}+\frac{3\lambda}{4}}.
    \end{equation}
\end{proposition}
\begin{proof}
    We use interpolation again to derive this proposition. First directly using Proposition 3.6 and \eqref{2.6}, for $|I|\le N-19$, we have
    \begin{equation}
        |ZZ^Iu(t,x)|+|SZ^Iu(t,x)|\lesssim \langle t+r\rangle^{-\frac{3}{4}} (C^3 \epsilon^{\frac{1}{8}} \langle \tau \rangle ^{\frac{5}{8}+p}+10)\lesssim C^3  \langle t+r\rangle^{-\frac{1}{8}+p}.
    \end{equation}
    Next using Proposition 3.4, for $0\le\lambda\le\frac{1}{100}$, we have 
    \begin{equation}
    \begin{split}
         |ZZ^Iu(t,x)|+|SZ^Iu(t,x)|& \le  \langle t+r\rangle |\partial Z^Iu(t,x)| 
         \\ &\lesssim C^{3}\epsilon^{\lambda^2} {\langle \tau+r\rangle}^{\frac{3\lambda}{2}}{\langle \tau-r\rangle}^{-\frac{1}{8}}
         \\& \le C^{3}\epsilon^{\lambda^2} {\langle \tau+r\rangle}^{\frac{3\lambda}{2}}.
    \end{split}
    \end{equation}
    Interpolating the two inequalities above, we derive that 
    \begin{equation}
      |ZZ^Iu(t,x)|+|SZ^Iu(t,x)|\lesssim   C^3\epsilon^{\frac{\lambda^2}{2}}{\langle \tau+r\rangle}^{-\frac{1}{16}+\frac{3\lambda}{4}+\frac{p}{2}}.
    \end{equation}
    Finally, using Lemma 2.3, we have
    \begin{equation}
        |G_{a}Z^Iu| \lesssim \langle t + |x| \rangle^{-1} ( |SZ^Iu| + | Z Z^Iu| )\lesssim C^3\epsilon^{\frac{\lambda^2}{2}}{\langle \tau+r\rangle}^{-\frac{17}{16}+\frac{3\lambda}{4}+\frac{p}{2}}.
    \end{equation}
\end{proof}

\subsection{End of the proof of Theorem 1.1} \label{sec 3.4}
In this section, we complete the proof of Theorem 1.1, starting with the proof of Proposition 3.1.
\begin{proof}
    \textbf{Step 1.} Close conformal energy estimate of $u$ and $\partial u$ in \ref{3.1}.
    \\For $|I|\le N-10$, act $ Z^I$ on both sides of wave equation in \eqref{add 1.2}, using Lemma 2.1, we have
    \begin{equation}
        -\Box Z^I u=\sum_{|I_1|+|I_2|\le |I|} C_{I_1,I_2,I}Q_0( Z^{I_1} u, Z^{I_2}w),
    \end{equation}
where $C_{I_1,I_2,I}$ are constants.
  Applying the  conformal energy estimate, Lemma 2.1, Lemma 2.2, and Claim 3, we have
      \begin{equation}
      \begin{split}
          \mathcal{F}^{\frac{1}{2}}(t,  Z^Iu) & \lesssim \mathcal{F}^{\frac{1}{2}}(0,   Z^Iu) + \sum_{|I_1|+|I_2|\le |I|}\int_{0}^{t} \left\| \langle \tau + r \rangle Q_0( Z^{I_1} u , Z^{I_2}w) \right\|_{L^{2}}  d\tau
          \\ &+ \sum_{a,b}\sum_{|I_1|+|I_2|\le |I|}\int_{0}^{t} \left\| \langle \tau + r \rangle Q_{ab}( Z^{I_1} u , Z^{I_2}w) \right\|_{L^{2}}  d\tau
          \\ & \lesssim \frac{200}{\epsilon}+ \sum_{|I_1|+|I_2|\le N-10}\int_{0}^{t} \left\| (|SZ^{I_1}u|+|ZZ^{I_1}u|) |ZZ^{I_2} w |\right\|_{L^{2}}  d\tau.
      \end{split}  \label{3.38}
\end{equation}
Note $|I_1|\le N-10$, using Hardy's inequality (Corollary 2.5), we have
\begin{equation}
\begin{split}
     &  \left\|\frac{S Z^{I_1}u}{r^{\frac{3}{10}}} \right\|_{L^2}\lesssim  \left\|{S Z^{I_1}u} \right\|_{L^2}^{\frac{7}{10}} \left\|{\partial S Z^{I_1}u} \right\|_{L^2}^{\frac{3}{10}} \le C \epsilon ^{-\frac{7}{10}} \langle \tau \rangle ^{2p},
     \\ &   \left\|\frac{Z Z^{I_1}u}{r^{\frac{3}{10}}} \right\|_{L^2}\lesssim  \left\|{Z Z^{I_1}u} \right\|_{L^2}^{\frac{7}{10}} \left\|{\partial Z Z^{I_1}u} \right\|_{L^2}^{\frac{3}{10}} \le C \epsilon ^{-\frac{7}{10}}\langle \tau \rangle ^{2p}.
\end{split}
\end{equation}
Note $|I_2|\le N-10$, using Proposition 3.5, we have 
\begin{equation}
      \left\|{Z Z^{I_1}w}{r^{\frac{3}{10}}} \right\|_{L^\infty}\lesssim   C^4 \epsilon^{\frac{5}{7}}\langle \tau  \rangle^{-\frac{6}{5}} \label{add 3.445}.
\end{equation}
Therefore, we have 
\begin{equation}
\begin{split}
    &  \sum_{|I_1|+|I_2|\le N-10}\int_{0}^{t} \left\| (|SZ^{I_1}u|+|ZZ^{I_1}u|) |ZZ^{I_2} w |\right\|_{L^{2}}  d\tau 
    \\ & \lesssim  \sum_{|I_1|+|I_2|\le N-10}\int_{0}^t \left\|{Z Z^{I_1}w}{r^{\frac{3}{10}}} \right\|_{L^\infty}\left\|\frac{S Z^{I_1}u}{r^{\frac{3}{10}}} \right\|_{L^2} d \tau\\ & + \int_{0}^t\left\|{Z Z^{I_1}w}{r^{\frac{3}{10}}} \right\|_{L^\infty} \left\|\frac{Z Z^{I_1}u}{r^{\frac{3}{10}}} \right\|_{L^2} d\tau
    \\ & \lesssim \int_0 ^t C^4 \epsilon^{\frac{5}{7}}\langle \tau +r \rangle^{-\frac{6}{5}+2p} C \epsilon ^{-\frac{7}{10}}
    \\ & \lesssim C^5 \epsilon ^{\frac{1}{70}}.
\end{split}
\end{equation}
 Substituting into \eqref{3.38}, we have 
 \begin{equation}
      \mathcal{F}^{\frac{1}{2}}(t,  Z^Iu)\lesssim \frac{1}{\epsilon}+ C^6 \epsilon ^{\frac{1}{70}}. \label{3.42}
 \end{equation}
 We may choose $\epsilon$ small enough to close the assumption of conformal energy on $Z^I u$.
 
 To close the conformal energy estimate for $\partial Z^Iu$, $|I|\le N-10$, we do almost the same thing, act $\partial Z^Iu$ on the wave equation of \eqref{add 1.2}, we obtain
    \begin{equation}
     -\Box \partial Z^I u=\sum_{|I_1|+|I_2|\le |I|}C_{I_1,I_2,I}Q_0( \partial Z^{I_1} u ,Z^{I_2}w)
     + E_{I_1,I_2,I}Q_0( Z^{I_1} u , \partial Z^{I_2}w),    
    \end{equation}
    where $C,D$ are constants.
 Applying the conformal energy estimate, Lemma 2.1, Lemma 2.2, and Claim 3, we have
      \begin{equation}
      \begin{split}
        &  \mathcal{F}^{\frac{1}{2}}(t,  \partial Z^Iu) 
        \\&\lesssim \mathcal{F}^{\frac{1}{2}}(0,   \partial Z^Iu) + \sum_{|I_1|+|I_2|\le |I|}\int_{0}^{t} \left\| \langle \tau + r \rangle Q_0( \partial Z^{I_1} u , Z^{I_2}w) \right\|_{L^{2}}  d\tau
          \\ &+\int_{0}^{t} \left\| \langle \tau + r \rangle Q_0( Z^{I_1} u , \partial Z^{I_2}w) \right\|_{L^{2}}  d\tau
          \\ & \lesssim 10+ \sum_{|I_1|+|I_2|\le N-10}\int_{0}^{t} \left\| (|S\partial Z^{I_1}u|+|Z\partial Z^{I_1}u|) |ZZ^{I_2} w |\right\|_{L^{2}}  d\tau
          \\ & +  \sum_{|I_1|+|I_2|\le N-10}\int_{0}^{t} \left\| (|S Z^{I_1}u|+|Z Z^{I_1}u|) |\partial ZZ^{I_2} w |\right\|_{L^{2}}  d\tau.
      \end{split}  \label{3.44}
\end{equation}
Note that the last line is totally the same as the term emerging in closing conformal energy of $Z^I u$(repeating the same treatment from \eqref{3.38} to \eqref{3.42}, note that the corresponding \eqref{add 3.445} is still valid since we only require $|I|\le N-8$ in Proposition 3.5 ), therefore, we have 
 \begin{equation}
     \sum_{|I_1|+|I_2|\le N-10}\int_{0}^{t} \left\| (|S Z^{I_1}u|+|Z Z^{I_1}u|) |\partial ZZ^{I_2} w |\right\|_{L_{x}^{2}}  d\tau\lesssim C^5 \epsilon ^{\frac{1}{70}}.\label{3.45}
 \end{equation}
 It only remains to treat the first term, apply Proposition 3.5, note $|I_1|\le N-10$, $|I_2|\le N-10$, and we have 
 \begin{equation}
     \begin{split}
       &  \sum_{|I_1|+|I_2|\le N-10}\int_{0}^{t} \left\| (|S\partial Z^{I_1}u|+|Z\partial Z^{I_1}u|) |ZZ^{I_2} w |\right\|_{L^{2}}  d\tau
    \\ & \lesssim        \sum_{|I_1|+|I_2|\le N-10} \int_{0}^{t} \|
    S\partial Z^{I_1}u\|_{L^{2}} \|ZZ^{I_2} w \|_{L^{\infty}} + \|Z\partial Z^{I_1}u\|_{L^{2}} \|ZZ^{I_2} w \|_{L^{\infty}} d\tau
    \\ & \lesssim  \sum_{|I_1|+|I_2|\le N-10} \int_{0}^{t} C \langle\tau\rangle ^p C^4 \epsilon^{\frac{5}{7}}\langle \tau  \rangle^{-\frac{3}{2}} d\tau
    \\ & \lesssim C^5 \epsilon^{\frac{5}{7}} \label{3.46}.
z     \end{split}
 \end{equation}
Substituting \eqref{3.45}, \eqref{3.46} into \eqref{3.44}, we derive that
\begin{equation}
     \mathcal{F}(t,  \partial Z^Iu)^{\frac{1}{2}}\lesssim C^5 \epsilon ^{\frac{1}{70}}+ C^5 \epsilon^{\frac{5}{7}}+10.\label{3.47}
\end{equation}
\textbf{Step 2}: Close the energy estimate of $Z^I u$ and $\partial Z^Iu$.

We start with $\partial Z^I u$, $|I|\le N-1$.
\\ Acting $\partial Z^I$ on the wave equation in \eqref{add 1.2}, by Lemma 2.7,  
  \begin{equation}
        -\Box \partial Z^I u=\sum_{|I_1|+|I_2|\le |I|}C_{I_1,I_2,I}Q_0( \partial Z^{I_1} u , Z^{I_2}w)+D_{I_1,I_2,I}Q_0( Z^{I_1} u , \partial Z^{I_2}w),
    \end{equation}
where $C,D$ are constants. By applying Lemma 2.7, Lemma 2.2, and Claim 2, we have 
\begin{equation}
\begin{split}
     &  \mathcal{E}^{\frac{1}{2}}(t,\partial Z^I u)
     \\ &\lesssim   \|\partial\partial Z^I u\|_{L^2}(0)+  \sum_{|I_1|+|I_2|\le N-1}  \int_0^t \|Q_0( \partial Z^{I_1} u , Z^{I_2}w)\|_{L^2} 
       \\ & +\|Q_0( Z^{I_1} u , \partial Z^{I_2}w)\|_{L^2} d\tau
    \\ &  \lesssim 2 \epsilon +   \sum_{|I_1|+|I_2|\le N-1} \int_0^t \|\partial \partial Z^{I_1} u \partial Z^{I_2}w\|_{L^2}+\|\partial Z^{I_1} u \partial \partial Z^{I_2}w\|_{L^2} d\tau \label{3.49}.
\end{split}
\end{equation}
For the first term. we have 
\begin{equation}
    \begin{split}
     & \sum_{|I_1|+|I_2|\le N-1}  \int_0^t \|\partial \partial Z^{I_1} u \partial Z^{I_2}w\|_{L^2}d\tau 
     \\ &\lesssim  \sum_{\substack{|I_1|\le N-1 \\ |I_2|\le N-15}} \int_0^ t\|\partial \partial Z^{I_1} u \partial Z^{I_2}w\|_{L^2} d\tau
        +\sum_{\substack{|I_1|\le N-18 \\ |I_2|\le N-1}} \int_0^ t\|\partial \partial Z^{I_1} u  \partial Z^{I_2}w\|_{L^2} d\tau \label{3.50}.
    \end{split}
\end{equation}
Using \eqref{3.1} and Proposition 3.3, for $|I_1|\le N-1$, $|I_2|\le N-15$, we have 
\begin{equation}
\begin{split}
   &  \sum_{\substack{|I_1|\le N-1 \\ |I_2|\le N-15}} \int_0^ t\|\partial \partial Z^{I_1} u \partial Z^{I_2}w\|_{L^2} d\tau  
      \\  & \lesssim   \sum_{\substack{|I_1|\le N-1 \\ |I_2|\le N-15}}\int_0^ t\|\partial \partial Z^{I_1} u \|_{L^2}\|\partial Z^{I_2}w\|_{L^\infty} d\tau  
      \\ & \lesssim C^3  \epsilon^2 \int_0^t \langle \tau\rangle ^{-\frac{3}{2}+p}d\tau
      \\ & \lesssim C^3 \epsilon^2 \label{3.51}.
\end{split}
\end{equation}
Using \eqref{3.1} and taking $\lambda=\frac{p}{3}$ in Proposition 3.4, for $|I_1|\le N-18$, $|I_2|\le N-1$, we have 
\begin{equation}
\begin{split}
   &  \sum_{\substack{|I_1|\le N-18 \\ |I_2|\le N-1}} \int_0^ t\|\partial \partial Z^{I_1} u \partial Z^{I_2}w\|_{L^2} d\tau  
      \\  & \lesssim   \sum_{\substack{|I_1|\le N-18 \\ |I_2|\le N-1}} \int_0^ t\|\partial \partial Z^{I_1} u \|_{L^\infty}\|\partial Z^{I_2}w\|_{L^2} d\tau  
      \\ & \lesssim C^{4}  \epsilon^{1+\frac{p^2}{9}} \int_0^t \langle \tau\rangle ^{-1+\frac{p}{2}}d\tau
      \\ & \lesssim C^{4}  \epsilon^{1+\frac{p^2}{9}} \langle t\rangle ^{\frac{p}{2}}\label{3.52}.
\end{split}
\end{equation}
Substituting \eqref{3.51}, \eqref{3.52} into \eqref{3.50}, we obtain that
\begin{equation}
     \sum_{|I_1|+|I_2|\le N-1}  \int_0^t \|\partial \partial Z^{I_1} u \partial Z^{I_2}w\|_{L^2}d\tau \lesssim C^{4}  \epsilon^{1+\frac{p^2}{9}} \langle t \rangle ^{\frac{p}{2}}+ C^3 \epsilon^2 \label{add 3.57}.
\end{equation}
For the second term, we have 
\begin{equation}
\begin{split}
   & \sum_{|I_1|+|I_2|\le N-1} \int_0^t \|\partial Z^{I_1} u \partial \partial Z^{I_2}w\|_{L^2} d\tau \\ &\lesssim \sum_{\substack{|I_1|\le N-1\\ |I_2|\le N-16}}\int_0^t \|\partial Z^{I_1} u \partial \partial Z^{I_2}w\|_{L^2} d\tau+\sum_{\substack{|I_1|\le N-17\\ |I_2|\le N-1}}\int_0^t \|\partial Z^{I_1} u \partial \partial Z^{I_2}w\|_{L^2} d\tau \label{3.54}.
\end{split}
\end{equation}
Using Corollary 2.6 for $\partial Z^{I_1}u$, and Proposition 3.4, for $|I_1|\le N-1$, $|I_2|\le N-16$, we have
\begin{equation}
    \begin{split}
     &   \sum_{\substack{|I_1|\le N-1\\ |I_2|\le N-16}}\int_0^t \|\partial Z^{I_1} u \partial \partial Z^{I_2}w\|_{L^2} d\tau 
     \\ & \lesssim  \sum_{\substack{|I_1|\le N-1\\ |I_2|\le N-16}}\int_0^t \|\frac{\partial Z^{I_1} u }{r^\frac{1}{10}}\|_{L^2}\|\partial \partial Z^{I_2}w r^{\frac{1}{10}}\|_{L^\infty} d\tau
     \\ &\lesssim  \sum_{\substack{|I_1|\le N-1\\ |I_2|\le N-16}}\int_0^t \|\partial Z^{I_1} u\|_{L^2}^{\frac{9}{10}}\|\partial\partial Z^{I_1} u\|_{L^2}^{\frac{1}{10}}\|\partial \partial Z^{I_2}w r^{\frac{1}{10}}\|_{L^\infty} d\tau
     \\ & \lesssim C^3\epsilon^{1+\frac{1}{10}} \int_0^t \langle \tau\rangle^{-\frac{3}{2}+\frac{1}{10}}d\tau
     \\ & \lesssim C^3 \epsilon^{\frac{11}{10}} \label{3.55}.
    \end{split}
\end{equation}
Using \eqref{3.1}, Proposition 3.4, where we take $\lambda=\frac{p}{3}$, for $|I_1|\le N-17$, $|I_2|\le N-1$, we have 
\begin{equation}
    \begin{split}
      &  \sum_{\substack{|I_1|\le N-17\\ |I_2|\le N-1}}\int_0^t \|\partial Z^{I_1} u \partial \partial Z^{I_2}w\|_{L^2} d\tau 
      \\ & \lesssim \sum_{\substack{|I_1|\le N-17\\ |I_2|\le N-1}} \int_0^t \|\partial Z^{I_1} u \|_{L^{\infty}}\|\partial \partial Z^{I_2}w\|_{L^2} d\tau 
      \\ & \lesssim C^{4-\frac{2p^2}{9}}\epsilon^{1+\frac{p^2}{9}} \int_0^t\langle \tau \rangle^{-1+\frac{p}{2}} d\tau
      \\ & \lesssim C^{4-\frac{2p^2}{9}}\epsilon^{1+\frac{p^2}{9}} \langle t\rangle ^{\frac{p}{2}} \label{3.56}.
    \end{split}
\end{equation}
Substituting \eqref{3.55}, \eqref{3.56} into \eqref{3.54}, we obtain that
\begin{equation}
     \sum_{|I_1|+|I_2|\le N-1} \int_0^t \|\partial Z^{I_1} u \partial \partial Z^{I_2}w\|_{L^2} d\tau \lesssim C^3 \epsilon^{\frac{11}{10}}+C^{4-\frac{2p^2}{9}}\epsilon^{1+\frac{p^2}{9}} \langle t\rangle ^{\frac{p}{2}} \label{3.57}.
\end{equation}
Combining \eqref{add 3.57}, \eqref{3.57}, using \eqref{3.49}, we have
\begin{equation}
\begin{split}
     \mathcal{E}^{\frac{1}{2}}(t,\partial Z^I u) & \lesssim  10\epsilon +C^3 \epsilon^{\frac{11}{10}}+C^{4-\frac{2p^2}{9}}\epsilon^{1+\frac{p^2}{9}} \langle t\rangle ^{\frac{p}{2}}+C^{4-\frac{p^2}{9}}  \epsilon^{1+\frac{p^2}{9}} \langle t \rangle ^{\frac{p}{2}}+ C^3 \epsilon^2
     \\ & \lesssim10\epsilon+ C^4 \epsilon^{1+\frac{p^2}{9}} \langle t\rangle ^{\frac{p}{2}}. \label{final step 2}
\end{split}
\end{equation}
By repeating the almost same process, we can also strictly improve the bound of $\mathcal{E}^{\frac{1}{2}}(t, Z^I u)$.
\\ \textbf{Step 3} Close the energy  estimate of $w$.
\\ For$|I|\le N$, acting $Z^I$ on both sides of Klein-Gordon part of \eqref{add 1.2}, we have 
\begin{equation}
    -\Box Z^I w +Z^I w = \sum_{|I_1|+|I_2|\le |I|}  C_{I_1,I_2,I}Q_0(Z^{I_1}u, Z^{I_2}w),
\end{equation}
where $C_{I_1,I_2,I}$ are constants. By applying the ghost weight estimate(Lemma 2.8), Claim 4, Lemma 2.2, we have 
\begin{equation}
    \begin{split}
        E_{gst}(t,Z^Iw)&\lesssim  \int_{\mathbb{R}^3} (|\partial_t Z^Iw|^2 + \sum_a |\partial_a Z^Iw|^2 + m^2 |Z^Iw|^2) \, dx(0)\\& + \sum_{|I_1|+|I_2|\le |I|} \int_0^t \int_{\mathbb{R}^3} | Q_0(Z^{I_1}u, Z^{I_2}w)\partial_t Z^Iw | \, dxd\tau
       \\& \lesssim 10\epsilon+\sum_a\sum_{|I_1|+|I_2|\le N} \int_0^t \int_{\mathbb{R}^3}|G_aZ^{I_1}u| |\partial Z^{I_2}w||\partial_t Z^Iw | \, dxd\tau
       \\& + \sum_a\sum_{|I_1|+|I_2|\le N} \int_0^t \int_{\mathbb{R}^3}  |\partial Z^{I_1}u| |G_a Z^{I_2}w||\partial_t Z^Iw | \, dxd\tau.
    \end{split} \label{3.60}
\end{equation}
When the good derivatives hit $Z^{I_1}u$, we have 
\begin{equation}
    \begin{split}
      &  \sum_a\sum_{|I_1|+|I_2|\le N} \int_0^t \int_{\mathbb{R}^3}|G_aZ^{I_1}u| |\partial Z^{I_2}w||\partial_t Z^Iw | \, dxd\tau
      \\ & \lesssim \sum_a\sum_{\substack{|I_1|\le N \\ |I_2|\le N-15}} \int_0^t \int_{\mathbb{R}^3}|G_aZ^{I_1}u| |\partial Z^{I_2}w||\partial_t Z^Iw | \, dxd\tau 
       \\ & + \sum_a\sum_{\substack{|I_1|\le N-19 \\ |I_2|\le N}} \int_0^t \int_{\mathbb{R}^3}|G_aZ^{I_1}u| |\partial Z^{I_2}w||\partial_t Z^Iw | \, dxd\tau.
    \end{split}   \label{3.61}
\end{equation}
Using Proposition 3.3, \ref{3.1}, Corollary 2.6, we have
\begin{equation}
\begin{split}
    & \sum_a\sum_{\substack{|I_1|\le N \\ |I_2|\le N-15}}\int_0^t \int_{\mathbb{R}^3}|G_aZ^{I_1}u| |\partial Z^{I_2}w||\partial_t Z^Iw | \, dxd\tau
    \\ & \lesssim \sum_a\sum_{\substack{|I_1|\le N\\ |I_2|\le N-15}} \int_0^t   \| \frac{\partial Z^{I_1}u}{r^{\frac{1}{4}}}\|_{L^2}\|r^{\frac{1}{4}} \partial Z^{I_2}w\|_{L^\infty} \|\partial_t Z^{I}w\|_{L^2}(\tau)d\tau
    \\ & \lesssim \sum_a\sum_{\substack{|I_1|\le N\\ |I_2|\le N-15}} \int_0^t   \|{\partial \partial Z^{I_1}u}\|_{L^2}^{\frac{1}{4}}\|{\partial Z^{I_1}u}\|_{L^2}^{\frac{3}{4}}\|r^{\frac{1}{4}} \partial Z^{I_2}w\|_{L^\infty} \|\partial_t Z^{I}w\|_{L^2}(\tau)d\tau
    \\ & \lesssim C^{4}\epsilon^{\frac{9}{4}} \int_0 ^t\langle \tau \rangle^{-\frac{5}{4}+p} d\tau
    \lesssim C^4 \epsilon^{\frac{9}{4} }. \label{3.62}
\end{split}
\end{equation}
Choosing $\lambda=\frac{1}{100}$ in Proposition 3.7, using \eqref{3.1}, we also have
\begin{equation}
\begin{split}
  &  \sum_a\sum_{\substack{|I_1|\le N-19 \\ |I_2|\le N}} \int_0^t \int_{\mathbb{R}^3}|G_aZ^{I_1}u| |\partial Z^{I_2}w||\partial_t Z^Iw | \, dxd\tau
  \\ & \lesssim \sum_a\sum_{\substack{|I_1|\le N-19 \\ |I_2|\le N}} \int_0^t \int_{\mathbb{R}^3}\|G_aZ^{I_1}u\|_{L^\infty} \|\partial Z^{I_2}w\|_{L^2}\|\partial_t Z^Iw \|_{L^2} d\tau
  \\& \lesssim  C^5 \epsilon^{\frac{\lambda^2}{2}+2}\int_0^t\langle \tau\rangle^{-\frac{17}{16}+\frac{p}{2}+\frac{3\lambda}{4}}d\tau
  \\&\lesssim C^5 \epsilon^{\frac{\lambda^2}{2}+2}.
\end{split}
     \label{3.63}
\end{equation}
Substituting \eqref{3.62} and \eqref{3.63} into \eqref{3.61},we obtain
\begin{equation}
       \sum_a\sum_{|I_1|+|I_2|\le N} \int_0^t \int_{\mathbb{R}^3}|G_aZ^{I_1}u| |\partial Z^{I_2}w||\partial_t Z^Iw | \, dxd\tau \lesssim C^5 \epsilon^{\frac{\lambda^2}{2}+2}+C^4 \epsilon^{\frac{9}{4} } \label{3.83}.
\end{equation}
When the good derivatives hit $Z^{I_2}w$, we have 
\begin{equation}
\begin{split}
    \sum_a\sum_{|I_1|+|I_2|\le N} \int_0^t \int_{\mathbb{R}^3}  |\partial Z^{I_1}u| |G_a Z^{I_2}w||\partial_t Z^Iw | \, dxd\tau
     \\ \lesssim \sum_a\sum_{\substack{|I_1|\le N \\ |I_2|\le N-15}} \int_0^t \int_{\mathbb{R}^3}  |\partial Z^{I_1}u| |G_a Z^{I_2}w||\partial_t Z^Iw | \, dxd\tau
     \\+\sum_a\sum_{\substack{|I_1|\le N-17 \\ |I_2|\le N}} \int_0^t \int_{\mathbb{R}^3}  |\partial Z^{I_1}u| |G_a Z^{I_2}w||\partial_t Z^Iw | \ dxd\tau. \label{3.84}
\end{split}    
\end{equation}
 For $|I_1|= N$, $|I_2|\le N-15$, on the one hand, using Proposition 3.3, \eqref{3.1}, Corollary 2.6, we have
\begin{equation}
    \begin{split}
       & \int_{\mathbb{R}^3}  |\partial Z^{I_1}u| |G_a Z^{I_2}w||\partial_t Z^Iw | \, dxd\tau
       \\ & \lesssim \sum_{\substack{|I_1|\le N-1 \\ |I_2|\le N-15}}  \int_{\mathbb{R}^3}  |\partial\partial Z^{I_1}u \langle \tau +r\rangle| |G_a Z^{I_2}w||\partial_t Z^Iw | \, dxd\tau
       \\ &\lesssim \sum_{\substack{|I_1|\le N-1 \\ |I_2|\le N-15}}  \|{\partial \partial Z^{I_1}u}\|_{L^2} \|\partial_a Z^{I_2 }w\langle \tau +r \rangle\|_{L^\infty}\|\partial_t Z^Iw \|_{L^2} \
        \\ & \lesssim C^4 \epsilon^3 \langle \tau \rangle^{-\frac{1}{2}+p}.
    \end{split}
\end{equation}
On the other hand, using Proposition 3.3, 
\begin{equation}
    \begin{split}
       & \int_{\mathbb{R}^3}  |\partial Z^{I_1}u| |G_a Z^{I_2}w||\partial_t Z^Iw | \, dxd\tau
       \\ &\lesssim \sum_{\substack{|I_1|\le N \\ |I_2|\le N-15}}  \|{ \partial Z^{I_1}u}\|_{L^2} \|\partial_a Z^{I_2 }w\|_{L^\infty}\|\partial_t Z^Iw \|_{L^2} \
        \\ & \lesssim C^4 \epsilon^2 \langle \tau \rangle^{-\frac{3}{2}+p}.
    \end{split}
\end{equation}
Interpolating the above two inequalities, 
\begin{equation}
  \sum_{\substack{|I_1|=N\\ |I_2|\le N-15}} \int_{\mathbb{R}^3}  |\partial Z^{I_1}u| |G_a Z^{I_2}w||\partial_t Z^Iw | \, dxd\tau  \lesssim C^4 \epsilon^{\frac{9}{4}} \langle \tau \rangle ^{-\frac{5}{4}+p}.\label{add 3.91}
\end{equation}
When $|I_1|\le N-1$, using Hardy's inequality, we have, 
\begin{equation}
    \begin{split}
       & \sum_a\sum_{\substack{|I_1|\le N-1 \\ |I_2|\le N-15}} \int_0^t \int_{\mathbb{R}^3}  |\partial Z^{I_1}u| |G_a Z^{I_2}w||\partial_t Z^Iw | \, dxd\tau
       \\ &\lesssim \sum_a\sum_{\substack{|I_1|\le N-1 \\ |I_2|\le N-15}} \int_0^t  \|\frac{\partial Z^{I_1}u}{r^{\frac{1}{4}}}\|_{L^2} \|\partial_a Z^{I_2 }wr^{\frac{1}{4}}\|_{L^\infty2}\|\partial_t Z^Iw \|_{L^2} \,d\tau
        \\ &\lesssim \sum_a\sum_{\substack{|I_1|\le N-1 \\ |I_2|\le N-15}} \int_0^t  \|\partial\partial Z^{I_1}u\|_{L^2}^{\frac{1}{4}} \|\partial Z^{I_1}u\|_{L^2}^{\frac{3}{4}}\|\partial_a Z^{I_2 }wr^{\frac{1}{4}}\|_{L^\infty2}\|\partial_t Z^Iw \|_{L^2} \,d\tau
       \\ & \lesssim C\epsilon^{\frac{9}{4}}\int_0^t \langle\tau\rangle^{-\frac{5}{4}+p}d\tau
       \\ & \lesssim  C\epsilon^{\frac{9}{4}}. \label{3.85*}
    \end{split}
\end{equation}
Therefore, we have
\begin{equation}
    \begin{split}
       & \sum_a\sum_{\substack{|I_1|\le N \\ |I_2|\le N-15}} \int_0^t \int_{\mathbb{R}^3}  |\partial Z^{I_1}u| |G_a Z^{I_2}w||\partial_t Z^Iw | \, dxd\tau
       \\ &\lesssim \sum_a \big(\sum_{\substack{|I_1|\le N-1 \\ |I_2|\le N-15}} + \sum_{\substack{|I_1|= N \\ |I_2|\le N-15}} \big)\int_0^t  |{\partial Z^{I_1}u}| |\partial_a Z^{I_2 }w||\partial_t Z^Iw | \,d\tau
        \\ & \lesssim \sum_a \int_0^t C^4 \epsilon^{\frac{9}{4}} \langle \tau \rangle ^{-\frac{5}{4}+p} d\tau + C\epsilon^{\frac{9}{4}}
        \\ &\lesssim C^4\epsilon^{\frac{9}{4}}.
        \label{3.85}
    \end{split}
\end{equation}
Choosing $\lambda=\frac{1}{100}$ in Proposition 3.4, using \eqref{3.1} and Cauchy's inequality, we have
\begin{equation}
    \begin{split}
      & \sum_a\sum_{\substack{|I_1|\le N-17 \\ |I_2|\le N}} \int_0^t \int_{\mathbb{R}^3}  |\partial Z^{I_1}u| |G_a Z^{I_2}w||\partial_t Z^Iw | \ dxd\tau 
      \\ & \le\sum_a\sum_{\substack{|I_1|\le N-17 \\ |I_2|\le N}} \int_0^t \int_{\mathbb{R}^3} C^{3}\epsilon^{\lambda ^2} \langle \tau +r\rangle^{-\frac{1}{2}-\frac{1}{16}+\frac{3\lambda}{2}}  \frac{|G_a Z^{I_2}w|}{\langle \tau -r\rangle^{\frac{1}{2}+\frac{1}{16}}}|\partial_t Z^Iw | \ dxd\tau
      \\ & \lesssim \sum_a\sum_{\substack{|I_1|\le N-17 \\ |I_2|\le N}} \int_0^t C^{3}\epsilon^{\lambda ^2} \langle \tau \rangle^{-\frac{1}{2}-\frac{1}{16}+\frac{3\lambda}{2}} \| \frac{G_a Z^{I_2}w}{\langle \tau -r\rangle^{\frac{1}{2}+\frac{1}{16}}}\|_{L^2}\|\partial_t Z^Iw \|_{L^2}d\tau
      \\ &\lesssim C^4 \epsilon^{1+\lambda ^2} \sum_a (\int_0^t\langle \tau \rangle^{-{1}-\frac{1}{8}+3\lambda}d\tau)^{\frac{1}{2}}(\int_0^t\| \frac{G_a Z^{I_2}w}{\langle \tau -r\rangle^{\frac{1}{2}+\frac{1}{16}}}\|_{L^2}^2)^{\frac{1}{2}}
      \\ &\lesssim C^5 \epsilon^{2+\lambda ^2} . \label{3.86}
    \end{split}
\end{equation}
Substituting \eqref{3.85}, \eqref{3.86} into \eqref{3.84}, we have
\begin{equation}
    \sum_a\sum_{|I_1|+|I_2|\le N} \int_0^t \int_{\mathbb{R}^3}  |\partial Z^{I_1}u| |G_a Z^{I_2}w||\partial_t Z^Iw | \, dxd\tau\lesssim C^5 \epsilon^{2+\lambda ^2}+C\epsilon^{\frac{9}{4}}.\label{3.87}
\end{equation}
By \eqref{3.83}, \eqref{3.87}, and \eqref{3.60}, we finally derive 
\begin{equation}
    E_{gst}(t,Z^Iw)\lesssim C^5 \epsilon^{2+\lambda ^2}+C\epsilon^{\frac{9}{4}}+C^5 \epsilon^{\frac{\lambda^2}{2}+2}+C^4 \epsilon^{\frac{9}{4} }+\epsilon^2+10\epsilon. \label{3.88}
\end{equation}
where $\lambda=\frac{1}{100}$.
By \eqref{3.42}, \eqref{3.47}, 
\eqref{final step 2}and \eqref{3.88} we have strictly improved all the estimates in \eqref{3.1}.Therefore for all data $(u_0,u_1,w_0,w_1)$ satisfying \ref{1.3}, $T_*=\infty$, thus the proof of Theorem 1.1 completes.

\end{proof}

\newpage
\bibliography{ref}
\end{document}